
\input amstex
\define\Date{April 19, 2000}
\documentstyle{amsppt}
\magnification=1200
\loadeurm
\loadeusm
\loadbold
\font\twbf=cmbx10 scaled\magstep2
\font\chbf=cmbx10 scaled\magstep1
\font\hdbf=cmbx10


\pagewidth{6.5 true in}
\pageheight{8.9 true in}
\nologo
\overfullrule=3pt
\define\0{^{@,@,\circ}}
\define\1{^{-1}}
\define\8{^{\infty}}

\define\A{I^{[q]} \Col cx^q}

\define\Ai{I^{[q_i]} \Col cx^{q_i}}
\define\Ais{(I^{[q_i]})^* \Col cx^{q_i}}

\define\Ass{\text{Ass}\,}
\define\Ap{I^{[q']} \Col cx^{q'}}

\define\Bigskip{\quad\bigskip\bigskip\bigskip}
\define\blank{\underline{\phantom{J}}}

\define\cF{{\Cal F}}

\define\clo#1#2#3{{{#1}^{#2}}_{#3}}

\define\cln{\colon}
\define\col#1{\,\colon_{#1}}
\define\Col{\,\colon\,}

\define\cQ{{\Cal Q}}

\define\cont{\supseteq}

\define\demb#1{\demo{\bf #1}}

\define\e{\epsilon}
\define\F{\boldkey F}
\define\I{I^{[q]}}

\define\imp{\Rightarrow}
\define\inc{\subseteq}

\define\Iq{I^{[q]}}

\define\Let{Let $R$ be a reduced Noetherian ring of positive
prime characteristic $p$}

\define\Min{\text{Min}\,}
\define\N{{\Bbb N}}
\define\notinc{\nsubseteq}
\define\pa{\parindent 20 pt}
\define\part#1{\item{\hbox{\rm (#1)}}}
\define\pl{^\boldkey{+}}
\define\PP{{\Bbb P}}
\define\q{^{[q]}}
\define\Q{^{[Q]}}

\define\sat{^{\hbox{sat}}}
\define\sc{;\,}

\define\T{(\Iq)^* \Col cx^q}
\define\vect#1#2{#1_1,\, \ldots, #1_{#2}}
\define\Z{{\Bbb Z}}
\quad\bigskip\bigskip
\centerline{\twbf LOCALIZATION AND TEST EXPONENTS}
\vskip 10 pt plus .5 pt minus .5 pt
\centerline{\twbf FOR TIGHT CLOSURE}

\topmatter
\author
by Melvin Hochster and Craig Huneke
\endauthor
\rightheadtext{LOCALIZATION AND TEST EXPONENTS FOR TIGHT CLOSURE}
\leftheadtext{Melvin Hochster and Craig Huneke}
\thanks{The authors were supported in part by grants from the National
Science Foundation.} \endthanks
\thanks{The second author wishes to thank the Max Planck Institute in
Bonn and the Fulbright Foundation for partial support of the
research for this manuscript.}\endthanks
\thanks{Version of \Date.}\endthanks
\endtopmatter

\document
{\baselineskip = 16 pt
\quad\bigskip
{\chbf
\centerline{\hdbf 1. INTRODUCTION}
}
\bigskip
   We introduce the notion of a {\it test exponent} for tight closure,
and explore its relationship with the problem of showing that
tight closure commutes with localization, a longstanding open question.
Roughly speaking, test exponents exist if and only if tight
closure commutes with localization:  mild conditions on the ring are
needed to prove this.     We give other, independent, conditions
that are necessary and sufficient for tight closure to commute
with localization in the general case, in terms of behavior of certain
associated primes and behavior of exponents needed to annihilate
local cohomology.  While certain related conditions (the ones given
here are weaker)
were previously known to be sufficient, these are the first
conditions of this type that are actually equivalent.

The difficult
calculation of \S4 uses associativity of multiplicities and many
other tools to show that sufficient conditions for localization
to commute with tight closure can be given in which asymptotic
statements about lengths of modules defined using the iterates
of the Frobenius endomorphism replace the finiteness conditions
on sets of primes introduced in \S3.
The result is local and requires special conditions on the rings:  one is
that countable prime avoidance holds.  This is not a very
restrictive condition however:  it suffices, for example,  for
the ring to contain an uncountable field.  Countable prime avoidance
also holds in any complete local ring.  But we also need the
existence of a {\it strong} test ideal (see the beginning of \S4).
We expect that, in the long run, this condition will also turn out not
to be very restrictive:  a strong test ideal for a reduced ring is known
to exist if every irreducible component of Spec$\,R$ has a resolution of
singularities obtained by blowing up an ideal that defines the
singular locus, and it is expected that this is always true in
the excellent case. Moreover, by very recent results, strong test ideals
always exist for complete reduced local rings.

We note that the reader may find other results related to
localization of tight closure in [AHH], [Hu3], [Ktzm1-2], [LySm], and
[Vrac1-2].

\Bigskip
{\chbf
\centerline{\hdbf 2. TEST EXPONENTS AND ASSOCIATED PRIMES}
}
\bigskip
\demb{Discussion 2.1: basic terminology and notation}
We shall assume throughout that $R$ is a Noetherian ring of positive
prime characteristic $p$, although this hypothesis is usually
repeated in theorems and definitions.  Moreover, because tight closure
problems are unaffected by killing nilpotents, we shall assume, unless
otherwise specified, that $R$ is reduced.  We shall usually assume
that the reduced ring $R$ has a test element $c$ (see the discussion
below).   We recall
some terminology and notation.  $R\0$ denotes the complement of
the union of the minimal primes of $R$, and so, if $R$ is reduced,
$R\0$ is simply the multiplicative system of all nonzerodivisors in $R$.
We shall write $\F^e$ (or $\F^e_R$ if we need to specify the base ring)
for the {\it Peskine-Szpiro} or {\it Frobenius} functor from $R$-modules
to $R$-modules.  This is a special case of the base change functor
from $R$-modules to $S$-modules that is simply given by
$S \otimes_R\,\blank$:
in the case of $\F^e$, the ring $S$ is $R$, but the map $R \to R$
that is used for the algebra structure is the $e\,$th iteration $F^e$
of the Frobenius endomorphism:  $F^e(r) = r^{p^e}$.
Thus, if $M$ is given as the
cokernel of the map represented by a matrix $\bigl(r_{ij}\bigr)$,
then $\F^e(M)$ is the cokernel of the map represented by the
matrix $\bigl(r_{ij}^{p^e}\bigr)$.  Unless otherwise indicated, $q$
denotes $p^e$ where $e \in \N$.  For $q = p^e$,
$\F^e(R/I) \cong R/I\q$,  where $I\q$ denotes the ideal generated
by the $q\,$th powers of all elements of $I$ (equivalently, of
generators  of $I$).  Note that $\F^e$ preserves both freeness and
finite generation of modules, and is exact precisely when
$R$ is regular (cf. [Her], [Kunz]). If $N \inc M$  we write $N\q$ for
the image of $F^e(N)$ in $F^e(M)$, although it depends on the inclusion
$N \to M$,  not just on $N$.  If $u \in M$ we write $u^{p^e}$
for the image $1 \otimes u$ of $u$ in $F^e(M)$. With this
notation,  $(u+v)^q = u^q + v^q$ and $(ru)^q = r^qu^q$ for
$u,\,v \in M$ and $r \in R$.

It is worth noting that, for any multiplicative system
$W$ in $R$,  $\F^e_{R_W}(M_W) \cong F^e_R(M)_W$,  where
$\blank\,\,_W$ indicates localization with respect to $W$.
In fact, for any  $R$-algebra $S$, one has that
$$
\F^e_S(S \otimes_R M) \cong S \otimes_R\F^e_R(M).
$$ Furthermore, $N\q_W \inc \F^e(M)_W$ may be
canonically identified with $(N_W)\q$.

An element $u \in M$, where $M$ is a finitely generated
$R$-module, is in the tight closure $N^*$  of $N \inc M$  if
there exists $c \in R\0$ such that $cu^q \in N\q$  for
all $q = p^e \gg 0$.   By the right exactness of tensor,
$\F^e(M/N) \cong \F^e(M)/N\q$,  where $q = p^e$.  It follows
easily that an element of $M$ is in $N^*$ if and only if
its image in $M/N$ is in $0^*$ in $M/N$.  It follows that in
considering whether  $u \in M$ is in $N^*$,  we may replace
$M$ by a finitely generated free module $G$ mapping onto $M$,  $N$
by its inverse image in $G$, and $u$ by any element
of $G$ that maps to $u$.

An element $c \in R\0$ is called {\it a test element} if, whenever $M$ is
a finitely generated $R$-module and $N \inc M$ is a submodule, then
$u \in M$ is in the tight closure of $N$ if and only if for
all $q = p^e$,  $cu^q \in N\q$ (the image of $\F^e(N) \to \F^e(M))$.
Thus, if the ring has a test element, it ``works" in any tight closure
test where some choice of $c \in R\0$ ``works."  Test elements
are also characterized as the elements of $R\0$ that annihilate
$N^*/N$ for all submodules $N$ of all finitely generated modules $M$.

A test element is called {\it locally stable} if its image in every local
ring of $R$ is a test element (this implies that it is a locally stable
test element in every localization of $R$ at any multiplicative
system).  A test element is called {\it completely stable} if
its image in the completion of each local ring of $R$ is
a test element:  a completely stable test element is easily
seen to be locally stable.  We refer the reader to
[HH1, \S6 and \S8 ],
[HH2], [HH5, \S6], and [AHH, \S2] for more information
about test elements and to \S3 of [AHH] for a discussion of
several basic issues related to the localization problem for
tight closure.

We note that, by some rather hard theorems, test elements are known
to exist.  For example, if $R$ is any reduced ring essentially of
finite type over an excellent local ring, then $R$ has a test element.
In fact, if  $c$ is any element of $R\0$ such that $R_c$ is regular
(and such elements always exist if $R$ is excellent and reduced),
then  $c$ has a power that is a completely stable test element.
This follows from Theorem (6.1a) of [HH5], and we shall make
use of this freely throughout.

Another important property of tight closure in characteristic $p$ is
that if  $R \to S$ is a ring homomorphism and $u \in N^*$ in $M$,
then the image $1\otimes u$ of $u$ in $S \otimes_R M$ is in the
tight closure, over $S$,  of
$$
\hbox{\rm Im}\,(S \otimes_RN \to S\otimes_R M)
$$
under very mild assumptions. This phenomenon
is referred to as {\it the persistence of tight closure}.
In particular, by
Theorem (6.23) (which, with the same proof, is valid whenever
$R$ is essentially of finite type over an excellent local ring
--- the result is stated only for the case of finite type) and
Theorem (6.24) of [HH5], if  $R$ is essentially of finite type over an
excellent local ring, if $S$ has a completely stable test element,
or if $R\0$ maps into $S\0$ (e.g., if  $R \to S$ is an inclusion
of domains or is flat)  then one has persistence of tight closure
for the ring homomorphism  $R \to S$.
\enddemo

\demb{Definition 2.2} Let $R$ be a reduced Noetherian ring of positive
prime characteristic $p$.  Let $c$ be a fixed test element for
$R$.  Let $N \inc M$ be a pair of finitely generated $R$-modules.
We shall say that  $q = p^e$  is a {\it test exponent} for $c,\, N, M$
if whenever $cu^Q \in N^{[Q]}$ and $Q \geq q$,  then $u \in N^*$.
In case $N$ is an ideal it is usually assumed that $M = R$,  and in
that case we speak of a test exponent for  $c, I$, with $R$ understood
to be the ambient module for  $I$. \enddemo

It is not all clear whether to expect test exponents to exist.
In this paper we shall prove, roughly speaking, that test exponents
exist, in general, if and only if tight closure commutes with
localization.  The question of whether tight closure commutes with
localization is open in general, but it is known in many important
special cases, and thus our results imply that test exponents do exist
rather often.  We expect this notion to be of great importance
because, whenever one can compute what the test exponent is, one
obtains an effective test for tight closure.

We are also hopeful that focusing attention on the problem of the
existence of test exponents may lead to a solution of
the localization problem for tight closure.  We should want to point
out that, if tight closure commutes with localization, then it
commutes with arbitrary smooth base change: for a precise statement,
see Theorem (7.18a) of [HH5].

To demonstrate one connection, we prove the easier half of the
result at once --- this half is implicit in [McD].

\proclaim{Proposition 2.3} \Let\ and let $N \inc M$ be finitely
generated $R$-modules.  Let $c$ be a locally stable test element for $R$.
If  there is a test exponent for  $c,\, N,\, M$ then for every
multiplicative system $W$ of $R$,  $N_W^*$ in $M_W$ over $R_W$
is $(N^*)_W$ (i.e., tight closure for the pair $N \inc M$ commutes
with localization).\endproclaim

\demb{Proof} The only problem is to show that if $u \in M$
and $u/1 \in (N_W)^*$ then  $u \in (N^*)_W$ (any element of $M_W$, after
multiplication by a suitable unit, is in the image of $M$).
Let  $q$ be a test exponent for $c,\,N,\,M$.  Then we can choose
$cu^q/1 \in N\q_W$ and we can choose  $f \in W$  such that
$fcu^q \in N^q$, and so $c(fu)^q \in N\q$.  But then $fu \in N^*$,  and
so  $u \in (N^*)_W$. \qed
\medskip
What is much less obvious is that a converse holds.

\proclaim{Theorem 2.4} \Let\ and
let $N \inc M$ be a pair of finitely generated $R$-modules.
Suppose that for every prime $Q \in \Ass (M/N^*)$,
$(N_Q)^* = (N^*)_Q$.  Then for every test element $c \in R\0$ such
that  $c$ is a test element in each of the rings $R_Q$ for every
associated prime $Q$ of $M/N$,
there is a test exponent for $c,\,M,\,N$. In particular,
if   $(N_Q)^* = (N^*)_Q$ for all associated primes $Q$ of
$M/N^*$ and $c$ is a locally stable test element for $R$,  then
$c, N, M$ has a test exponent.

Thus, if tight closure
commutes with localization for the pair $N \inc M$ at
associated primes of $M/N^*$, and $R$ has a locally stable test element,
then tight closure  commutes with localization in general for the pair $N
\inc M$. \endproclaim

We postpone the proof until we have established some preliminary
results that make the argument transparent (it is given
immediately following the proof of Theorem 2.7).  First:

\proclaim{Proposition 2.5} \Let\ and let $M$ be a finitely generated
$R$-module.  Let $N$, $N'$,  and $N_i$,  where $i$ varies in an index set,
be submodules of $M$.  Let $c \in R\0$ be a test element.
{\pa
\part{a} $q = p^e$ is a test exponent for $c,\,N,\,M$ if and only
if it is a test exponent for $c,\,0,\,M/N$.
\part{b} If $q = p^e$ is a test exponent for $c,\,N,\,M$ then so is
every larger power of $p$.\par
\part{c} If $N \inc N' \inc N^*$ and $q$ is a test exponent for
$c,\,N',\,M$ then it is a test exponent for $c,\,N,\,M$.
In particular, this holds when $N' = N^*$.
\part{d} If $N \inc N' \inc N^*$,
$d$ is a test element, and $q$ is a test
exponent for $cd,\,N,\,M$, then
$q$ is a test exponent for $c, \,N',\,M$.
\part{e} If  $N_1$  and $N_2$
have the same tight closure in $M$,  and $W$ is a multiplicative
system consisting of test elements,  then  $c_1,\,N_1, \,M$ has a test
exponent for every choice of test element $c_1 \in W$  if and only if
$c_2,\, N_2, \,M$ has a test exponent for every choice of test element
$c_2 \in W$.
\part{f} If $q$ is a test exponent for $c, \,N_i, M$  for every
index $i$  (the index set may be infinite), every $N_i$ is tightly
closed, and  $N = \bigcap_i\,N_i$, then $q$ is a test exponent for
$c,\,N,\, M$.
\part{g} If $\vect N h$ are submodules of $N$ such that there is
a test exponent for each of $c,\,N_i, M$  for $1 \leq i \leq h$,
then there is a test exponent for $c,\, \bigcap_i N_i,\, M$.\par
}
\endproclaim
\demb{Proof}  Part (a) is immediate from the definitions of test
exponent and tight closure, while (b) is immediate from the definition
of test exponent. In part (c),  since  $N'$ and $N$ have the same
tight closure,  if $q$ is a test exponent for $c,\,N',\,M$ then
$cu^Q \in N\Q$ for $Q \geq q$ implies that $cu^Q \in {N'}\Q$,  and
so $u \in {N'}^* = N^*$.  For part (d),  if  $q$ is a test
exponent for $cd,\,N,\,M$,  $Q \geq q$ and $cu^Q \in {N'}\Q$,
then since  $N'$ is in the tight closure of $N$,  we have that
${N'}\Q$ is in the tight closure of $N\Q$, and so $cdu^Q \in N\Q$,
and we can conclude that $u \in N^* = {N'}^*$.  For part (e) we
might as well assume that $N_2$ is the tight closure of $N_1$:
it suffices to compare each with its tight closure.  The result
is immediate from the combination of parts (c) and (d).  Part
(f) is immediate from the definition, for if $cu^Q \in N\Q$
then $cu^Q \in N_i\Q$ for all $i$,  and so $u \in N_i$ for all $i$,
which shows that $u \in N$.  Part (g) follows from (b) and (f):
we may use the supremum of the finitely many test exponents for
the various $N_i$. \qed \enddemo

\proclaim{Proposition 2.6} \Let\ and let $N \inc M$ be finitely
generated modules such that $M/N$ has a unique associated
prime $P$, such that $N$ is tightly closed  in $M$,  and such
that $N_P$ is tightly closed in $M_P$ over $R_P$.  Let $c \in R\0$ be
any test element for $R$ that is also a test element in
$R_P$.  Then there is a test exponent for $c,\, N,\, M$. \endproclaim

\demb{Proof} We first consider the case where $P$ is maximal and
$R = R_P$.  Let $N_e$ denote
the set of elements  $u \in M$  such that  $cu^q \in (N\q)^F$, where
$q = p^e$ and $^\F$ indicates Frobenius closure in $M\q = \F^e(M)$.
Clearly,  $N^* \inc N_e$,
for every $e$. We claim that  $N_{e+1} \inc N_e$,  for if $cu^{pq} \inc
{N^{[pq]}}^F$,  then  $c^{q'}u^{pqq'} \in N^{[pqq']}$,  which
certainly implies that $c^{pq'}u^{pqq'} \in N^{[pqq']}$,
and this shows that $cu^q \in {N\q}^F$,  as required.  Thus,
this sequence of modules is eventually constant, since  $M/N$
has finite length in this case.  But once  $N_e = N_{e+1} = \ldots$,
the common value must be $N^*$,  for if  $u$  is in all of these,
we have that  $cu^Q \in {N^{[Q]}}^F \inc {N^{[Q]}}^*$ for
all $Q \geq p^e$,  and so  $c^2u^Q \in N^{[Q]}$ for all $Q \geq p^e$,
as required.  But then $cu^{q'} \in N^{[q']}$ for $q' \geq p^e$ implies
that $u \in N_e = N^*$, as required.

In the general case, note that the hypothesis is stable when we localize
at $P$, and so, by the case already proved, there exists $q'$ such
that  $cu^{q'}/1 \in {N_P}^{[q']}$ implies that $u/1 \in N_P$.  But
then, since $cu^{q'} \in N^{[q']}$ is preserved when we localize at
$P$,  we have that this implies that $u$ is in the contraction of $N_P$ to
$M$,  and since  $M/N$ is $P$-coprimary, this implies that
$u \in N$. \qed \enddemo

\proclaim{Theorem 2.7 (existence of primary decompositions that respect
tight clos\-ure)}
\Let, and suppose that $R$ has a test element.
Let  $N \inc M$ be finitely generated modules and suppose
that  $R_P$  has a test element (e.g., is
excellent)  for every associated
prime of  $M/N$.  Suppose that $N$  is tightly closed in $M$ and
that  $N_P$  is
tightly closed in $M_P$ for every associated prime  $P$ of $M/N$.

Then $N$ has a primary decomposition in $M$ in which for every
associated prime $P$ of $M/N$,  the $P$-primary
component  is tightly closed, and remains tightly closed
(in $M_P$, over $R_P$) after localization at  $P$. \endproclaim

\demb{Proof}  It suffices to construct such a primary decomposition
for $N_P \inc M_P$ over $R_P$ for each associated prime  $P$ of $M/N$.
Take all these submodules of $M_P$
and contract them  to  $M$, with $P$ varying among the associated primes
of $M/N$. (Note that the contraction of a tightly
closed submodule of $M_P$ to $M$ is tightly closed in
$M$:  this is a consequence of the persistence of tight closure,
which is automatic for the flat homomorphism $R \to R_P$.)
When there are several primary
components for the same  prime  $Q$  (each tightly closed, each
remaining tightly closed over  $R_Q$),  intersect them all.

This gives  a primary decomposition of $N$ with the required properties.
To see that it is, in fact, a primary decomposition, call the intersection
$N'$.  If  $N'/N \not= 0$, it contains a nonzero element whose annihilator
is an associated prime of $M/N$, and this remains true after localizing at
that  associated prime, which gives a contradiction.

Thus,  there is no loss of generality in assuming that  $(R,P)$
is local, that   $P$  is an associated prime of $M/N$, and that we have
solved the problem of constructing suitable primary components
after localizing at any of the other associated primes,  by induction on
the dimension
of  $R_P$.    Thus,  we may give a primary decomposition over  $R_Q$ for
every associated prime  $Q$   strictly contained in  $P$,  and then intersect
the contractions of all the modules occurring as above.   Call the
intersection of these other primary components
$H$.  Thus, $H$ has a primary decomposition using modules that
are primary for the other associated primes of $M/N$, and
that are tightly closed and remain so upon localization at
the respective associated primes.

If we localize at any element of  $P$,   only the other associated
primes remain,  and so  $H/N$  is killed by a power of  $P$,  and
is a finite length module.

Now consider the descending chain of submodules  $(N + P^nM)^* = N_n$,
which are tightly closed submodules  of  $M$  containing  $N$.
Then  $(N_n \cap H)/N$  is contained in  $H/N$,   and so the chain
$N_n \cap H$  is eventually stable,  which means that
$N_n \cap H$  is eventually stable.  But the intersection of
the $N_n$   is  $N$,  because  $N$  is tightly closed and
$R$  has a test element (if  $u \in N_n$ for all $n$ and
$c$ is the test element,  then for all $q$,  $cu^q \in (N + P^nM)\q
\inc N\q + P^n\F^e(M)$, and since  $N\q$ is $P$-adically
closed in $\F^e(M)$,  we have that for all $q$,  $cu^q \in N\q$, so
that $u \in N^* = N$, as required).  The stable value of $N_n \cap H$
must be the same as
$\bigcap_n\,(N_n \cap H) = (\bigcap_n\, N_n) \cap H = N \cap H = N$.
Thus, for all $n \gg 0$,  $N_n \cap H = N$,  and we may use
$N_n$ for any sufficiently large $n$ as the required tightly
closed $P$-primary component. \qed \enddemo

\demb{Proof of Theorem 2.4} For the first statement,
by Theorem 2.7 coupled with Proposition 2.6,
$N$ is a finite intersection of submodules $N_i$ of $M$ such
that there is a test exponent for $c,\,N_i, M$,  and the result is
then immediate from Proposition 2.5, parts (c) and (g).

For the second statement choose a test element of $R$ and for each
associated prime $Q$ of $M/N^*$ an element of $R\0$ that maps to
a test element for $R_Q$.  The product will be a test element $c$ of
$R$ that is also a test element for every $R_Q$,  and so $c, N, M$
has a test exponent.
\qed\enddemo

\proclaim{Corollary 2.8} \Let\ and suppose that $R$ has a locally
stable test element.  Let $N \inc M$ be finitely generated $R$-modules
and
assume in addition that there is a submodule
$N'\inc N^*$ such that $N'_Q = N_Q$ at all
associated primes
$Q$ of $M/N^*$, and such that the tight closure of $N'$ does commute with
localization. Then the tight closure of $N$ commutes with localization.
\endproclaim

\demb{Proof} By Theorem 2.4, it suffices to prove that
$(N^*)_Q = (N_Q)^*$
for all primes $Q$ that are associated to $N^*$. Fix such a $Q$.
Then $(N_Q)^* = (N'_Q)^* = ({N'}^*)_Q\inc (N^*)_Q$. \qed\enddemo

\proclaim{Corollary 2.9} \Let\ and suppose that $R$ has a locally
stable test element.
Suppose that $I \inc R$ is generically a complete
intersection (i.e.,  there is an ideal $I'$ generated by
a regular sequence such that  $I_P = I'_P$ for every minimal
prime $P$ of $I$) and $I^*$ has no embedded primes.
Then $(I_W)^* = (I^*)_W$
for all multiplicatively closed sets $W$.
\endproclaim

\demb{Proof} This follows immediately from Corollary 2.8 and the fact that
localization commutes with tight closure for ideals generated by
regular sequences (see Theorem (4.5) of [HH5]). \qed\enddemo

\demb{Remark} Corollary 2.8 can be combined with numerous other
theorems on when tight closure commutes with localization to
give other results similar to Corollary 2.9.  For example, if
$R$ is a domain of acceptable type in the sense of [AHH, p.\ 87]
(a mild condition satisfied by homomorphic images of Cohen-Macaulay
rings and by algebras essentially of finite type over an excellent
local ring) then in Corollary (2.9) we only need the ideal
$I'$ to be of height $d$ and generated by $d$ elements:  cf.\
Theorem (8.3b) on p.\ 112 of [AHH]. (If $R$ is not a domain the result is
valid if one has the required hypotheses modulo every minimal prime.)

\proclaim{Corollary 2.10} \Let\ with a locally stable test element.
If $(J_Q)^* = (J^*)_Q$ for every primary ideal $J$ with radical $Q$,
then $(I_W)^* = (I^*)_W$
for every ideal $I$ such that $I^*$ has no embedded primes, and for
every multiplicatively closed set $W$.
\endproclaim

\demb{Proof}  First note that the hypothesis implies that
$(J_W)^* = (J^*)_W$ for every primary
ideal $J$ and for every multiplicative set $W$.

Fix an arbitrary ideal $I$ with $I^*$ unmixed. Let $\vect Q n$ be the
associated primes of $I^*$.  These are
the same as the minimal primes of $I^*$,  and since
$I$ and $I^*$ have the same radical, they are
precisely the minimal primes of $I$.  Let $J_i$ be the primary
component of $I$ corresponding to $Q_i$. Since the $Q_i$ are minimal over $I$,
this component is uniquely determined. By the Theorem 2.4,  it is
enough to prove that $(I^*)_{Q_i} = (I_{Q_i})^*$ for each $Q_i$.  We have that
$(I_{Q_i})^* = (J_i^*)_{Q_i}$. It suffices to prove that
$J_i^*\inc (I^*)_{Q_i}$, for then $(I_{Q_i})^* = (J_i^*)_{Q_i}
\inc (I^*)_{Q_i}$. Choose $u\in I:J_i$ such that $u\notin Q_i$.
Let $x\in J_i^*$. For some $d\in R\0$ and for all large $q$,
$dx^q\in (J_i)^{[q]}$. Multiplying by $u^q$ gives that
$d(ux)^q\in I^{[q]}$, proving that $ux\in I^*$ and $x\in (I^*)_{Q_i}$
as required. \qed
\enddemo

\demb{Discussion 2.11} \Let. Given the results of this section,
it is natural to ask how the associated primes
of $I \inc R$ and of $I^*$ are related, and, more generally, how
the associated primes of $M/N$ and $M/N^*$  are related when
$N \inc M$, a finitely generated $R$-module.
We want to comment that while the minimal primes are the
same,  there is, in general, no comparison, in either direction, for
the other associated primes.  The examples below show this
even in the case of ideals. (Note:  to see that a minimal prime $Q$ of
$M/N$ is still in the support of $M/N^*$ we may localize
at that prime,  since $(N_Q)^* \cont (N^*)_Q$, and then use the
fact that $(M/N)_Q$, if nonzero, maps on $R_Q/QR_Q$.)
Embedded primes of $I$ may fail to
be associated primes of $I^*$,  and embedded primes of $I^*$ may
fail to be associated primes of $I$.  In Example 2.12 below
$I$ has the maximal ideal as an embedded prime but $I^*$ does not.
In Example 2.13 below, $I$ is unmixed but $I^*$ has the maximal
ideal as an embedded prime.
\enddemo

\demb{Example 2.12} Let $R$ be a normal local ring of
positive prime characteristic and of dimension 3 that
is not Cohen-Macaulay (for definiteness, one may take the Segre
product of a homogeneous coordinate ring of an elliptic curve,
e.g., $K[x,\,y,\,z]/(x^3+y^3+z^3)$, with $K$ a field of
positive characteristic different from 3, with  $K[u,v]$, a homogeneous
coordinate ring for $\PP^1_K$, and localize at the irrelevant ideal,
i.e., at the unique maximal homogeneous ideal).  Let $f,\,g$ be
part of a system of parameters.  Since $R$ is not Cohen-Macaulay,
the maximal ideal is an embedded prime of $I = (f,\,g)R$.  By the
colon-capturing property of tight closure (see for example,
Theorem (1.7.4) of [HH6]), a third parameter for
$R$ is not a zerodivisor on $I^*$,  and so the maximal ideal of
$R$ is not an associated prime of $I^*$.  There are similar
examples in all dimensions. \enddemo

\demb{Example 2.13} Let $K$ be a field of positive prime characteristic
$p \not =3$.  Take
$$
R = K[X,\,Y,\,U,\,V]/(X^3Y^3 + U^3 + V^3) = K[x,\,y,\,u,\,v].
$$
Then $R$ is geometrically normal ($p \not= 3$) since the partial
derivatives of the defining polynomial  include $3U^2$, $3V^2$,
which form a regular sequence in the ring.
Let $I = (u,v, x^3)R$.  By the persistence of tight closure
(cf. the last paragraph of 2.1 above)
and the fact that $z^2$ is in $(u,\,v)^*$ in $K[z,\,u,\,v]/(z^3+u^3+v^3)$
(see, for example, the beginning of \S4 of [Ho1]),
we  have that $x^2y^2$ is in the tight closure of  $(u,v)$  and so
it is in $I^*$.  Now,  $R/I$ is isomorphic
with  $K[x,y]/(x^3)$ and so $I$ is unmixed,  while  $I^*$
contains  $x^2y^2$,  so that both $y^2$ and $x$  multiply
$x^2$ into $I^*$ (as well as  $u$ and $v$ of course).  This
will show that  $R/I^*$ has $m = (x,y,u,v)$ as an embedded
prime provided that  $x^2$ is not in $I^*$.  But this is
true even if we kill $u$, $v$, and $y$, and tight closure
persists under homomorphisms for affine algebras (cf.\ the last paragraph
of 2.1 again).
\enddemo

We next observe that if tight closure commutes
with localization after a faithfully flat extension, then it commutes
with localization.

\proclaim{Proposition 2.14} Let $R$ be a reduced Noetherian ring
of positive prime characteristic $p$ such that $R$ has a completely
stable test element.  Let $R \to S$ be faithfully flat.  Let $N \inc M$
be finitely generated $R$-modules.  Let $W$ be a multiplicative
system in $R$.  If tight closure commutes with localization for
the pair $S \otimes_R N \inc S \otimes_R M$ and the multiplicative
system that is the image of $W$ in $S$,  then it commutes with
localization for the pair $N \inc M$ and the multiplicative system $W$.
\endproclaim

\demb{Proof} It suffices to show that if $u \in M$  is such that
$u/1$ is in the tight closure of $N_W$ in $M_W$ then $u$ is in $(N^*)_W$.
But then $1 \otimes u$ is in the tight closure of $S_W \otimes_R N$ in
$ S_W \otimes_R M$ over  $S_W$,  by the persistence of tight closure
(this is trivial in the flat case),  and by our hypothesis this
implies that $1 \otimes u$ is in $(S \otimes_R N)^*$ in
$(S \otimes_R M)^*$ over $S$.   This implies that $u$ is in
$N^*$,  by Corollary (8.8) of [Ho2], p.\ 143. \qed \enddemo

\demb{Remark 2.15} The conclusion of (2.14) is valid under a substantial
weakening of the hypothesis on the homomorphism $R \to S$.  Instead
of being faithfully flat, it suffices if it preserves height in
the sense of condition $(*)$ of Corollary (8.8) of [Ho2] and
persistence of tight closure holds:  the
proof is the same.

\proclaim {Corollary 2.16} Let $R$ be a reduced ring and $K$ a field
such that $R$ satisfies at least
one of the following five conditions:
{\pa
\part{1} $R$ is finite type over $K$
\part{2} $R$ is essentially of finite type over $K$
\part{3} $R$ is an excellent local ring with residue field $K$
\part{4} $R$ is of finite type over an excellent local ring
with residue field $K$.
\part{5} $R$ is essentially of finite type over an excellent local
ring with residue field $K$.\par
}
Then $R$ has a faithfully flat extension $S$  that satisfies the
same condition, but such that $K$ is uncountable and $F$-finite,
and, in cases (3)--(5), such that the local ring is complete.

Hence, the question of whether tight closure commutes with
localization for any of the rings in the five classes discussed
above can be reduced to a corresponding case where the ring contains an
uncountable field, is $F$-finite, and, in cases (3), (5),
where the local ring is complete as well.  \endproclaim

\demo{Proof} The last statement follows from the next to last
statement by Proposition 2.14.  We may trivially replace the
local ring by its completion in (3)--(5).   The statements
about rings essentially of finite type follow from the corresponding
statements for rings of finite type.   In cases (3)--(5)
choose a coefficient field, and let $L$ be a field obtained by
adjoining uncountably many indeterminates to it (or simply to
$K$ in case (1)).  Replace the local ring by its complete tensor
product with $L$ in cases (3), (4).  One may then use the $\Gamma$
construction of \S6 of [HH5] to make an excellent, faithfully flat
extension of the local ring such that the residue field is $F$-finite,
without losing the property that its tensor product with $R$ over
the original local ring is reduced.  Finally, one may replace this
local ring by its completion.  Case (1) is simply the special case
where the local ring has dimension 0. \qed\enddemo

Note that there are several implications among the conditions listed
in Corollary 2.16 (e.g., (5) $\imp$ (4) $\imp$ (3) and (4) $\imp$ (2)
$\imp$ (1)):  we have stated the result as we did because it contains
the information that whichever one of these conditions holds, {\it that
particular condition} can be preserved as one modifies $R$ to contain an
uncountable field.

\demb{Remark 2.17} Let $R$ be any ring containing an
uncountable field $K$ and let $I$ be a finitely generated ideal of $R$.
Let $\{J_n: n \in \N\}$ be a countable family of ideals of $R$ whose
union contains $I$.  Then $I$ is contained in one of $I_n$.  To see
this, let $V$ be the finite-dimensional $K$-vector space spanned
over $K$ by a finite set of generators for $I$.  Then $V$ is covered
by the vector spaces $I_n$ as $n$-varies, and it suffices
to show that one of these is $V$, for if $V \inc I_n$ then $I \inc I_n$.
Thus, we have reduced to a question purely about vector spaces.
But the result for vector spaces is clear if the dimension of $V$ is at
most one,  and follows at once by induction from the fact that $V$
has uncountably many mutually distinct subspaces of codimension one,
each of which, by the induction hypothesis, will be contained in
at least one of the $I_n$. It follows that at least two of these
are contained in the same $I_n$,  and that forces $V \inc I_n$.

Thus, if a Noetherian ring $R$ contains an uncountable field,
and the ideal $I$ is not contained in any of the countably
many ideals $I_n$,  then $I$ has an element that is not in any
of the $I_n$.  If this property holds when the $I_n$ are prime,
we say that $R$ {\it has countable prime avoidance}.

This discussion, coupled with Corollary 2.16, shows
that, in the main cases, the question of whether tight closure
commutes with localization reduces to the case where the
ring has countable prime avoidance.  We assume countable
prime avoidance in the main result of \S4.

\demb{Remark 2.18} We also want to remark that every Noetherian
ring $R$ of positive prime characteristic $p$ has a faithfully
flat Noetherian extension containing an uncountable field $K$.
For any ring $R$, if $T$ is a set of indeterminates let
$R(T)$ denote the localization of the polynomial ring $R[T]$ in these
indeterminates at the multiplicative system of all polynomials
whose coefficients generate the unit ideal.  $R(T)$ is easily seen to be
faithfully flat over $R$.  Evidently,
if $R$ contains any field $K$, e.g., $\Z/pZ$,  then
$R(T)$ contains $K(T)$,  and so if $R$ has characteristic $p$ and
$T$ is uncountable then $R(T)$ contains the uncountable field $K(T)$.
But, for any Noetherian ring $R$,  $S = R(T)$ is Noetherian. (It is easy
to see that expansion and contraction gives a bijection between
maximal ideals of $R$ and those of $S$.
It suffices if every prime $P$ of $S$ is finitely generated.
Let $W$ run through the finite subsets of $T$ and let $P_W$ denote
the contraction of $P$ to $R(W) \inc R(T)$.  It is clear that
$R(W)$ is Noetherian and so each $P_W$ is a finitely generated
prime ideal of $R(W)$.  Now $P_WS$ is prime in $S$, since
$S$ is obtained from $R(W)$ by adjoining indeterminates and localizing,
and $P_WS \inc P_{W'}S$ if  $W \inc W'$.  The prime $P$ is the union
of the $P_WS$ as $W$ varies.  We claim that $P_W S = P$ for any
sufficiently large choice of $W$.  The point is that $P$ is contained
in a maximal ideal of $S$,  say  $mS$, where $m$ is a maximal ideal
of $R$  with, say,  $n$ generators. Then for all $W$,
$mS$ lies over $mR(W)$ which has height at most $n$.  From this
it follows that if one has a chain of the form
$P_{W_0}S \subset \, \cdots \, \subset P_{W_h}S$
in which the inclusions are strict, then $h \leq n$,  since the
inclusions will also be strict in
$P_{W_0}R(W) \subset \, \cdots \, \subset P_{W_h}R(W)$
for some sufficiently large $W \supseteq W_h$, and this contradicts
the fact that $mR(W)$ has height at most $n$.)
\enddemo

\demb{Remark 2.19} We note that countable prime avoidance holds
in any complete local ring by Lemma 3 of [Bur].
\enddemo
\vfill\eject
\Bigskip
{\chbf
\centerline{\hdbf
3. NECESSARY AND SUFFICIENT CONDITIONS} \smallskip
\centerline{\hdbf FOR LOCALIZATION}
}
\bigskip

In this section we give a necessary and sufficient condition for the tight
closure of an ideal to commute with localization. Several {\it sufficient}
conditions have
been given in previous papers. Most notable among these were
that tight closure
commutes with localization provided that the following two conditions hold
(cf.\ [Ktzm1-2]):
\roster
\item For every ideal $I$, the union, over $q$, of the sets of maximal
associated primes of $(\Iq)^*$ is a finite set.
\item For every prime ideal $P$, there exists a positive integer
$k$ such that for all $q = p^e$, $$P^{kq}H^0_P(R/((\Iq)^*)_P) = 0.$$
\endroster

However, these conditions are not known to be necessary. Theorem 3.5
below gives a pair of conditions that, together, characterize
precisely when tight closure commutes with localization for
an ideal  $I$.  The conditions like C2$^*$ discussed below
are clearly weaker than the second condition above. Our condition
$C1$ is reminiscent of the first condition above, in that it asserts the
finiteness of certain sets of primes, but it is not immediately
apparent how to compare it to (1) directly.

Several preliminaries are needed before stating our main result,
Theorem 3.5 below.  In particular,  we need to identify certain sets of
primes in $R$:
these are characterized by several equivalent conditions in the next
result.

\proclaim{Proposition 3.1} \Let\ with at least one locally stable
test element.
Let $I \inc R $ and $x \in R$ be fixed.
The following conditions are equivalent for a
prime ideal $Q$:
{\roster
\item  For some (equivalently, every) locally stable test element
$c \in R\0$, $Q$ is
minimal over infinitely many of the ideals $\Iq \Col cx^q$.
\item For some (equivalently, every) locally stable test element $c \in R\0$,
$Q$  is minimal over  $\Iq \Col cx^q$ for all $q \gg 0$.
\item For some (equivalently, every) locally stable test element $c \in R\0$,
$Q$ is minimal over infinitely many ideals $(\Iq)^* \Col cx^q$.
\item For some (equivalently, every) locally stable test element $c \in R\0$,
$Q$ is minimal over $(\Iq)^* \Col cx^q$ for all $q \gg 0$.
\item $x\notin (I_Q)^*$, and $x\in (I_P)^*$ for all $P\subsetneq Q$.
  \endroster}
\endproclaim

\demb{Proof} Fix a single locally stable test element $c$.
It will suffice to
show that the conditions in (1) -- (4) are equivalent to the
condition in (5) for that one choice of $c$, since (5) does not
refer to $c$.  Thus, throughout the rest of the proof, $c$ is
a fixed locally stable test element and when we refer to conditions
(1) -- (4)
we are referring to them for that single fixed test element $c$.

Clearly $(2) \imp (1)$ and $(4) \imp (3)$.

We next prove that $(*)$ if
$Q$ is minimal over $\Ai$ for all $i \geq 1$
then it is minimal over $\Ais$ for all
$i \gg 0$. Once we have shown this, it follows at once that $(2)
\imp (4)$ and $(1) \imp (3)$.
Hence once we establish $(*)$, we need
only show that $(5) \imp (2)$
and $(3) \imp (5)$ to complete the argument.

Suppose that $(*)$ is false. Since $\Ai\inc \Ais$,
$Q$ will be minimal over every $\Ais$ that it contains. Suppose
that $Q$ does not contain $\Ais$ for infinitely many $i$.
After localizing at $Q$, we have that
$cx^{q_i}\in ((I^{[q_i]})^*)_Q$, which implies that
$c^2x^{q_i}\in (I^{[q_i]})_Q$ for infinitely many $i$. This shows
that $x\in (I_Q)^*$
(cf.\ [HH2], Lemma (8.16), p.\ 79). But then for all large $q$,
$cx^q\in I_Q^{[q]}$, proving
$\A\nsubseteq Q$ for all large $q$.
This contradiction finishes the proof of
$(*)$.

Next we prove that $(3) \imp (5)$. Since $Q$ is minimal over infinitely
many of the ideals  $(\Iq)^* \Col cx^q$, $Q$ contains infinitely many of
the ideals $\A$, and hence $x\notin (I_Q)^*$. Let $P\subsetneq Q$. After
localizing at $P$, the assumption guarantees that $cx^q\in ((I^{[q]})^*)_P$
for infinitely many values of $q$, forcing $c^2x^q\in (I^{[q]})_P$ for the
same values of $q$.  As above  (cf.\ [HH2], Lemma (8.16), p.\ 79)
this implies that $x\in (I_P)^*$.

Finally assume (5). First observe that $Q$ must contain $\A$ for all
$q \gg 0$. If not, then after localizing at $Q$ we obtain that
$x\in (I_Q)^*$, a contradiction. Next, suppose that $P\subsetneq Q$.
By assumption, $x\in (I_P)^*$, which means that for all $q$, there
exist elements $w_q\notin P$ such that $w_qcx^q\in \I$. Then for all
$q$, $\A\nsubseteq P$. It follows that $Q$ is minimal over $\A$
whenever it contains $\A$. \qed\enddemo

\smallskip
\demb{Definition 3.2} \Let\ having at least one locally stable test element,
let $I \inc R$ be an ideal and let $x \in R$.
The primes satisfying the equivalent conditions  of (3.1)
we call
the {\it stable} primes associated to $I$ and $x$, and we denote
by $T_I(x)$ the set of stable primes associated to $I$ and $x$.
Let  $T_I = \bigcup_{x \in R} T_I(x)$.
\enddemo

\medskip

For a fixed ideal $I$ we consider the following two conditions:
{\roster
\item"C1:" For every $x\in R$, the set $T_I(x)$ is finite.
\item"C1$^*$:" The set $T_I$ is finite.
\endroster}

For a fixed ideal $I$ and locally stable test element $c$
we also consider the following two conditions:

{\roster
\item"C2:" For every $x\in R$, if $Q\in T_I(x)$, there exists an
integer $N$, possibly depending on $Q$, such that for all $q \gg 0$,
$Q^{Nq}\inc (\A)_Q$.
 \item"C2$^*$:" For every $x\in R$, if $Q\in T_I(x)$,
there exists an
integer $N$, possibly depending on $Q$, such that for all $q \gg 0$,
$Q^{Nq}\inc (\T)_Q$.
\endroster}
\medskip
We shall soon use these conditions (Theorem 3.5 below)
to characterize precisely
when tight closure commutes with localization for $I \inc R$.
However, we need some preliminary results, as well as some new
notation and terminology.

By a {\it square} locally stable test element $c$ we mean one
such that $c = d^2$ where $d$ is a locally stable test element.
For technical reasons that stem from Proposition 3.3(d) below,
it is often advantageous to work with a square locally stable
test element.

We denote by $\sqrt{J}$ the radical of the ideal $J$ and
by $\Min(J)$ the set of minimal primes of $J$.

\proclaim{Proposition 3.3} \Let\ with a locally stable test element,
let $I$ be an ideal of $R$, and let $x$ be an element of $R$.
{\pa
\part{a} A prime ideal $Q$ of $R$ has the property that
$x \notin (IR_Q)^*$ (in $R_Q$) if and only if $Q$ contains
an ideal in $T_I(x)$.  Hence,
if localization commutes with tight closure for $I \inc R$,
then $T_I(x)$ is the set of minimal primes of $I^*\Col Rx$ and, in
consequence, in this case, $T_I$ is the set of associated primes of $I^*$.
In particular, if $I \inc R$ is such that localization commutes with
tight closure then all of the sets $T_I(x)$ and even $T_I$ itself are
finite.
\part{b} $T_I(x) = T_{I^*}(x)$ and $T_I = T_{I^*}$.
\part{c} If $W$ is a multiplicative system in $R$ then
the elements of $T_{I_W}(x/1)$ (respectively, $T_{I_W}$), working
over $R_W$, are the expansions of those primes in
$T_I(x)$ (respectively, $T_I$) that do not meet $W$.
\part{d} If $c$ is a square locally stable test element then the
sequence $\sqrt{\A}$ is nonincreasing as $q$ increases.  Hence,
if $Q\in \Min (\Ap)$ and $\A\inc Q$ for some $q\leq q'$, then
$Q\in \Min (\A)$.
\part{e} If $c$ is a square locally stable test element and $P$
is a minimal prime of $I\q\Col cx^q$,  then $P$ contains an
element of $T_I(x)$.
\part{f} If $c$ is a square locally stable test element,
$T_I(x)$ is finite, and $q$ is so large that
all elements of $T_I(x)$ are minimal over $I\q\Col cx^q$ (which
is true for all $q \gg  0$),  then the minimal primes
of $I\q\Col cx^q$ are precisely the elements of $T_I(x)$, and
$\sqrt{I\q\Col cx^q} = \bigcap_{P \in T_I(x)}\,P$.
\part{g}
$T_I(x) = \emptyset$ iff $x \in I^*$.  Moreover,  if $c$ is a
test element then the condition that
$I\q\Col cx^q = R$ for all $q \gg 0$ (respectively, for all $q$)
is also equivalent.  \par
}
\endproclaim

\demb{Proof} (a) Given any prime $P$ such that $x \notin (IR_P)^*$,
every larger prime $Q$ has the same property, since there is
a flat homomorphism $R_Q \to R_P$ and tight closure is persistent.
Since $R$ has DCC on prime ideals, every such prime contains a minimal
such prime.  But the minimal primes with this property constitute
$T_I(x)$, by part (5) of Proposition 3.1. The other statements in
(a) are immediate.

Parts (b), and (c) follow at once from the
characterization of $T_I(x)$ in part (5) of Proposition 3.1.

For part (d),  let $c = d^2$ where $d$ is a locally stable
test element.  It suffices to prove that
$\sqrt{I^{[pq]}\Col cx^{pq}} \inc \sqrt {I\q\Col cx^q}$.
Let $\A\inc P$ and assume that $I^{[pq]}\Col cx^{pq}\notinc P$.
Then $d^2x^{pq}\in I^{[q']}_P$.
Hence $(dx^{q})^p\in (I^{[q]})^{[p]}_P$ and it follows that
$dx^{q}\in ((I^{[q]})_P)^*$. In particular, $d^2x^q\in (I^{[q]})_P$, which
means that $\A\notinc P$, contradiction. This proves the first assertion.

The second statement is immediate from the first: suppose that
$Q\in \text{Min}(\Ap)$ and $\A\inc Q$ for some $q\leq q'$. If $Q$ is
not minimal over $\A$, then  there is a prime $P$ such that
$\A\inc P\inc Q$, with $P\ne Q$. Then from the first part, $P$ must
also contain $\Ap$, and so $Q$ cannot be minimal over $\Ap$ either.

To prove (e), fix a minimal prime  $P = P_0$ of  $J_q = I\q \Col cx^q$.
Since $P \supseteq \sqrt{J_{pq}}$,  it contains a minimal prime
$P_1$ of $J_{pq}$.  Continuing in this way we get a sequence
$$
P_0 \supseteq P_1 \supseteq \,\cdots \supseteq P_i \supseteq \, \cdots
$$
such that for every $i$,  $P_i$ is a minimal prime of $J_{qp^i}$.
Since the prime ideals of a Noetherian ring have DCC,  we can choose
$n$ such that $P_i = P_n$ for all $i \geq n$,  and it follows that
$P_n$ is in $T_I(x)$ and is contained in $P$.

Part (f) is immediate, since,  once all of the finitely
many primes in $T_I(x)$ are minimal primes of $I\q\Col cx^q$,
part (e) shows that there cannot be any others.

Finally, to prove (g), note that, by part (f),
if $T_I(x) = \emptyset$, and one considers
$J_q = I\q\Col cx^q$ using a square locally stable test element $c$,
then the set of minimal primes of $J_q$ is empty for
all $q$.  This means that every $J_q = R$, which says that
$cx^q \in I\q$ for all $q$,  and this implies that
$x \in I^*$.  The other direction is clear.  The last two
conditions quite generally characterize when $x \in I^*$ when
$c$ is a test element.  \qed
\enddemo

\Let, let $I$ be an ideal of $R$, and let $x,\,y \in R$.  Let
$P$ be a prime ideal in $T_I(x)$.  We shall say that $y$ {\it
clears $P$ from $T_I(x)$} if $P$ is not in $T_I(xy)$.  We
shall say that the ideal $J \inc R$ clears $P$ from
$T_I(x)$ if every element of $J$ clears $P$ from $T_I(x)$.

The following result gives some basic facts about clearing.
Part (g) is a bit different, although analogous, and will
be needed in \S4.

\proclaim{Lemma 3.4 (clearing lemma)} \Let, let $I$ be an ideal of $R$,
and let $x,\,y \in R$, and let $c$ be a square locally stable test element
for $R$.
{\pa
\part{a} $y$ clears $P \in T_I(x)$ if and only if
for all $q \gg 0$,  $y^q/1 \in (I\q \Col cx^q)_P$.  Hence, the set
of all elements of $R$ that clear $P$ from $T_I(x)$ is an ideal
contained in $P$.
\part{b} Condition C2 or C2$^*$ for $I$ and $c$ implies that every prime
$P$ of $T_I(x)$ has a power that clears $P$ from $T_I(x)$.
\part{c} Suppose that $\vect P n$ are finitely many primes of
$T_I(x)$ such that each $P_i$ has a power that clears $P_i$ from
$T_I(x)$.  Let $\cQ$ be a family of primes $Q_j$ none of which
contains any of the $P_i$.  If the family $\cQ$ is finite
there is an element $y$  of $\cap_N P_i^N$ for large $N$ not in any
of the $Q_j$ that clears all of the $P_i$ from $T_I(x)$.
\part{d} With hypothesis as in (c) except that the family
$\cQ$ is countable, the same conclusion holds provided that
$R$ has countable prime avoidance.
\part{e} Suppose that $T_I(x)$ and $T_I(xy)$ are finite.  Then
$\bigcap_{P \in T_I(x)}\,P \inc \bigcap_{P \in T_I(xy)}\,P$,
and, if $y$ clears at least one element of $T_I(x)$,
the inclusion is strict.
\part{f} Every prime in $T_I(xy)$ contains a prime in $T_I(x)$.
\part{g} Let $(R, m)$ be local and let $z \in R$ be such that
$z^qH^0_m(R/I\q) = 0$ for all $q$.  Then for all $u \in R$ and all
$q$,  $m$ is not an associated prime of $I\q\Col uz^q$. \par
}
\endproclaim
\demb{Proof}  (a) The condition that $P$ not be in $T_I(xy)$ is simply
that $P$ not be a minimal prime of $I\q\Col c(xy)^q$ for
$q \gg 0$,  and since $P$ is minimal over $I\q \Col cx^q$,  which
is smaller,  this is equivalent to saying that  $P$ does
not contain $I\q\Col c(xy)^q$ for $q \gg 0$,  i.e., that
$(I\q\Col cx^qy^q)_P = R_P$ for $q \gg 0$, which holds iff
$y^q/1 \in (I\q\Col cx^q)_P$ for $q \gg 0$.  The second statement
is then obvious.

Part (b) for condition C2 is immediate from the definition.
We may use C2$^*$ instead, because $T_{I^*}(x) = T_I(x)$ and
$T_{I^*}(xy) = T_I(xy)$. Parts (c) and
(d) are obvious.

The first statement in (e) follows from the fact that if $T_I(z)$
is finite then $\bigcap_{P \in T_I(z)}\,P = \sqrt{I\q\Col cz^q}$
for all sufficiently large $q$ by 3.3(f), taken with the
obvious inclusion $I\q:cx^q \inc I\q\Col cx^qy^q$. The final
statement is then obvious.

For part (f), note that any prime $P$ in $T_I(xy)$ contains
the radical of $I\q\Col c(xy)^q$ for some $q$, and that
this contains the radical of $I\q \Col cx^q$, so that
$P$ contains at least one minimal prime of $I\q \Col cx^q$,
and each of these contains an element of $T_I(x)$ by 3.3(e).

For part (g) suppose that $m$ is associated to
$I^{[q]} \Col uz^q$, and choose $y\notin I^{[q]} \Col uz^q$ such that
$my\inc I^{[q]} \Col uz^q$, i.e., $mz^qyu \inc I\q$. Since $z^q$ is not
in any associated prime of $I\q$ except $m$,  it follows that if we
take a primary decomposition of $I\q$,  $myu$ is in the intersection $I'$
of the primary components for primes other than $m$. Since $I'/I\q$
is supported only at $m$,  $myu$ is killed by a power of $m$.
This shows that $yu$
represents an element in $H^0_m(R/I^{[q]})$,
and so $z^qyu\inc I^{[q]}$ i.e., $y \in I^{[q]} \Col uz^q$,
a contradiction.
\qed\enddemo

With these preliminary results our characterization becomes
straightforward:

\proclaim{Theorem 3.5} \Let\ with a square locally
stable test element $c \in R\0$, and let
$I\inc R$
be an ideal. The following are equivalent:
{\roster
\item  $(I_W)^* = (I^*)_W$ for all multiplicatively closed
sets $W$ in $R$.
\item Condition C1 (or C1$^*$) holds for $I$ and condition C2 (or C2*)
holds for $I$ and $c$.

\endroster}
\endproclaim

\demb{Proof} We have already observed that (1) implies
C1$^*$ which obviously implies C1:  see Proposition 3.3(a).
Thus, we may ignore
the parenthetical comment about C1$^*$ in the proof.

Assume (1) and suppose that
$Q\in T_I(x)$. Then $Q$ is a minimal prime of
$I^* \Col x$,  and so we can  chose $n \gg 0$ so that
$Q^n\inc (I^*\Col x)_Q$, i.e., $Q^n(x/1) \inc (I^*)_Q$.
Then for all $q$, $c(Q^n)^{[q]}x^q\inc I^{[q]}_Q$, i.e.,
$(Q^n)\q \inc (I\q \Col cx^q)_Q$.
If $t$ exceeds the number of generators of $Q^n$, then
$(Q^n)^{tq}\subseteq (Q^n)^{[q]}$, and so  C2 holds
with $N = nt$. Hence, (1)
$\imp$ C1 and C2.  Now C2$^*$ is obviously weaker than C2, so that
to complete the proof it suffices to show that C1 and C2$^*$ together
imply (1).

We assume otherwise and get a contradiction.
We may assume that $W = R - Q$
for some prime ideal $Q$ by the results of \S2 here (or using the result of
[AHH], Lemma (3.5), p.\ 79). Suppose
that $x\in (I_Q)^*$
but $x\notin (I^*)_Q$.
We may assume that $x\in R$.  By Noetherian induction,
among all $x$ giving a counterexample there is one such
that $\bigcap_{P \in T_I(x)} \, P$ is maximal.
By 3.3(g), since  $x \notin I^*$ this is not the unit ideal, i.e. $T_I(x)$
is not empty.  Let $P \in T_I(x)$ be any element.  Then
$P$ is not contained in $Q$, by 3.3(a), since $x \in (I_Q)^*$,
while using C2$^*$ we know from 3.4(b) that
$P^N$ clears $P$ from $T_I(x)$ for large $N$.  Thus,  we
may choose  $y \in P^N - Q$
such that $y$ clears $P$ from $T_I(x)$. Since $y$ is invertible in $R_Q$
we still have that $xy \notin I^*R_Q$, while
$xy \in (IR_Q)^*$ is clear.  But the intersection of
the primes in $T_I(xy)$ is strictly larger than the intersection
of those in $T_I(x)$, by 3.4(f),  and this contradicts the hypothesis for
the Noetherian induction.
\qed\enddemo

\demb{Remark 3.6} Suppose that $R$ is as in Theorem 3.5 and
that condition C2 holds in $R$ but that
localization fails to commute with tight closure for $I \inc R$.
The proof of Theorem 3.5 evidently shows that condition C1 fails
for $I$, $x$ such that $x \in (I_Q)^* - (I^*)_Q$ for some prime
$Q$ of $R$. \enddemo
\vfill\eject
\Bigskip
{\chbf
\centerline{\hdbf 4. GROWTH OF FROBENIUS IMAGES} \smallskip
\centerline{\hdbf AND LOCALIZATION}
}
\bigskip

In this section we focus our attention on the behavior of certain
functions related to the Hilbert-Kunz function that we believe
control the localization of tight closure in the local case.
We shall informally say that ``localization holds for $R$"
to mean that tight closure commutes with localization for
all ideals of $R$.

As observed at the end of \S2, to settle the localization problem
for, say, excellent reduced local rings, it suffices to handle the
complete case:  we may even assume that the ring is complete
with an uncountable residue field.   Because localization holds if
it holds modulo every minimal prime, it is sufficient to prove it
for the case of a complete local domain.

A complete local domain $R$ is module-finite
over a complete regular local ring $A$ with the same residue field.
For sufficiently large $q$ the extension of fraction fields corresponding
to the inclusion  $A^{1/q} \inc A^{1/q}[R]$ will be separable.
If we know that localization holds for $A^{1/q}[R]$ then it holds
for $R$, by Proposition 2.14 as extended in Remark 2.15.

Therefore, it is reasonable to study the problem of localization
of tight closure for a local domain $R$ module-finite over an excellent
regular local ring $A$,  and such that the extension of fraction fields is
separable:  we shall call such extensions {\it generically
\'etale}.  Since the case where the ring is complete implies
all of the most important local cases, little is lost by assuming
that countable prime avoidance holds, and likewise, as we shall
see in 4.2 below, we may as well assume that there is a strong test
ideal.

We shall attack the problem of localization in this situation
in the main result, Theorem 4.5, of this section.  Our goal is to show
in a ``minimal example" of the possible failure of tight closure
to commute with localization,
in the presence of a certain boundedness condition on
behavior of local cohomology, localization becomes equivalent
to an assertion about the asymptotic behavior of lengths of certain
sequences of
modules defined in terms of iterates of the Frobenius endomorphism.
The length conditions, surprisingly, replace the finiteness conditions
on the sets $T_I(x)$ discussed earlier.

We need several preliminaries.

\demb{Definition 4.1} \Let.  If $(R,m)$ is local we shall say
that condition (LC) holds for $R$ if for every ideal $J$,
there is an integer $N$ such that $m^{Nq}H^0_m(R/J\q) = 0$
for all $q$ (evidently, if one knows this for all $q \gg 0$,
it follows for all $q$ after enlarging $N$ if necessary).

If $I$, $x$ and a
locally stable test
element $c$ are fixed, for each $Q \in T_I(x)$
we define
$$
\lambda_q(Q) = \frac {\lambda(R_Q/(\A)_Q)} {q^{\text{dim}(R_Q)}},
$$
where $\lambda$ denotes length (in this case, over the local ring
$R_Q$).
\enddemo

\demb{Definition and discussion 4.2} \Let\ and let $J$ be an ideal of R. $J$
is called {\it a strong test ideal} for $R$ if
$J$ meets $R\0$,  and  $JI^* = JI$  for every ideal $I$ of $R$.
The main result of [Hu2] shows that there is such an ideal $J$, which
is also a defining ideal for the singular locus in Spec$\,R$, provided
that for every minimal prime $P$ of $R$,  the singularities of
$R/P$ can be resolved by blowing up an ideal that defines the
singular locus.  If blowing up such a $J_0$ resolves the singularities,
then any high power of $J_0$ will be a strong test ideal.  Such ideals
are not unique, although there is a largest one.

By a very recent result [Vrac1], if $R$ is a reduced local ring of
positive prime characteristic $p$ such that $R$ has a completely
stable test element, then if $(R,\,m)$ is complete or if the test ideal
$\tau$ is $m$-primary (and in many other cases),  then
$\tau$ is a strong test ideal!

Suppose that one has a strong test ideal $J$ with generators
$\vect j k$.
Then for every ideal $I$, if $u \in I^*$,  then  $uJ \inc IJ$,
and we get equations  $uj_s = \sum_{t=1}^k i_{s,t}j_t$.  The usual
determinantal trick for proving integral dependence on an ideal
shows that $u$ is integrally dependent on $I$ using an equation
of degree $k$,  since  $u$ is an eigenvalue of the $k \times k$
matrix $\bigl(\,i_{s,t}\,\bigr)$.  This gives an extremely
useful uniform bound on the degrees needed for equations displaying
that elements in tight closures of ideals are in the integral closure
of the ideal.

We refine the determinantal trick slightly.  Suppose that
$v$, $u$ are elements
of a domain $R$, and that $J \not= 0$ and  $I'$ are ideals such
that  $Ju \in J(Rv+I') = Jv + JI'$.  Then we obtain that
$u$ is an eigenvalue of a $k \times k$  matrix of the form
$\bigl(\,r_{s,t}v + i'_{s,t}\,\bigr)$ where the $r_{s,t} \in R$
and the $i'_{s,t} \in I'$.  Let $V$ be a new variable.
The characteristic
polynomial of $\bigl(uI - \bigl(\,r_{s,t}V + i'_{s,t}\,\bigr)\bigr)$
may be expanded as a polynomial of total
degree $k$ in the variables $U$ and $V$ that is monic in $U$:
the terms of degree smaller than
$k$ have coefficients in $I'$, which can be readily seen by thinking
modulo $I'$.  Thus, the characteristic polynomial yields a degree $k$
homogeneous polynomial
$P(U,\,V)$ with coefficients in $R$, monic in $U$,
such that $P(u,\,v) \in I'$.
\enddemo

Before proving our main theorem, two preliminary lemmas are needed.

\proclaim{Lemma 4.3} Let  $R$  be a reduced ring of prime characteristic
$p$, module-finite and torsion-free over
a regular
domain $A$. Let $\cF$ be the fraction
field of $A$ and suppose that $\cF\otimes_A R$ is  \'etale (i.e.,
separable) over $\cF$, i.e., that $R$ is generically \'etale
(or smooth) over $A$.   Then:
\pa
\part{a} For all $q$,  $A^{1/q}[R] \cong A^{1/q} \otimes_A R$ is flat
over $R$.  Moreover,
there exist elements $d \in A\0$ such that
 $dR\inc A[R^{p}]$ (equivalently, taking $p\,$th roots,
 $d^{1/p}R^{1/p} \inc A^{1/p}[R]$).

\part{b} If $d \in A\0$ is chosen so that $dR\inc A[R^{p}]$,
then for all $q \gg 0$,
$$d^{\frac {q-1} {q}}R^{1/q}\inc A^{1/q}[R].
$$
\part{c} If $c\in A\0$ is chosen so that $cR\inc A[R^p]$ (see (a) above)
then for every ideal $I$ of $R$ and all $q,\, q'$ we  have that
$$(I^{[qq']} \col R cx^{qq'})\cap A\inc (I^{[q']} \col R cx^{q'})^{[q]}.$$

\part{d} If $d \in A$ is such that $dR^{1/q} \inc A^{1/q}[R]$ then
for all $q'$,  $d^q(I^{[qq']}\Col x^{qq'}) \inc (I^{[q']}\Col x^{q'})\q$.

\part{e} Suppose, moreover, that  $R$ has a strong test ideal $J$ with
$k$ generators.  Suppose that $d \in A\0$ is such that
$dR \inc A[R^p]$. Let  $c = d^k$.  Let $c' \in R\0$.
Assume that there is an
integer $N > 0$ such that $c'm^{Nq}\inc (y^{q}) + (I^{[q]}\,\cln_Rx^{q})$
for all $q \gg 0$. Then there exists an integer $L > 0$ such that
for all large $q$, $m^{NLq}\inc (y^{q}) + (I^{[q]}\,\col R cx^{q})$.
\endproclaim

\demb{Proof} Both statements in part (a) follows from Lemma (6.4) on
p.\ 50 of [HH2]
and the discussion that immediately precedes it,  where it
is shown that if $c \in A\0$ is such that $R_c$ is \'etale over
$A_c$ then $c$ has a power $c_1$ such that $c_1R^{1/p} \inc A^{1/p}[R]$.
We may then take  $d = c_1^p$.

(b) Write $q = p^n$.
 Let $h_q =
(1 + p \,+\,\cdots\,+\, p^{n-1})/p^n$. We use induction on $n$ to prove that
$d^{h_q}R^{1/q}\inc A^{1/q}[R]$.  The case $n=1$ is immediate.
Assume the result for $n$ and take $p\,$th roots.
We get that
$d^{h_n/p}R^{1/qp}\inc A^{1/qp}[R^{1/p}]$. Multiplying by $d^{1/p}$
gives that
$$
d^{(l+h_n)/p}R^{1/qp}\inc A^{1/qp}[A^{1/p}[R]],
$$
(using the $n = 1$ case). Since $(l+h_q)/p = h_{pq}$, we have
completed the inductive step. Finally,
$1 + p \,+\,\cdots\,+\, p^{n-1} = (q-1)/(p-1) \leq q-1$, where $q = p^n$.

(c) Let $a\in I^{[qq']} \Col cx^{qq'}$ be in $A$ as well, so
that  $acx^{qq'} \in I^{[qq']}$ and
take $q\,$th roots to
obtain that
$$c^{1/q}a^{1/q}x^{q'}\in I^{[q']}R^{1/q}.$$
It follows from part (b) above that $c^{(q-1)/q}R^{1/q}\inc A^{1/q}[R]$ and
so $a^{1/q}cx^{q'}\in I^{[q']}A^{1/q}[R]$. From the flatness
of $A^{1/q}[R]$ over $R$ we then obtain that
$a^{1/q} \in I^{[q']} \col R cx^{q'}$  (since $a \in A$,
$a^{1/q} \in A^{1/q} \inc A^{1/q}[R]$), and now we may
take $q\,$th powers to get the stated result.

(d) If $ux^{qq'} \in I^{[qq']}$ then taking $q\,$th roots
yields that  $u^{1/q}x^{q'} \in I^{[q']}R^{1/q}$ and then multiplying
by $d$ shows that $du^{1/q}x^{q'} \in I^{[q']}A^{1/q}[R]$, and so with
$B = A^{1/q}[R]$,  we have
$du^{1/q} \in  I^{[q']}B\col B x^{q'}B$, and this is
$(I^{[q']}\col R x^{q'})B$ because $B$ is flat over $R$,
so that $du^{1/q} \in (I^{[q']}\col R x^{q'})A^{1/q}R$.  Taking
$q\,$th powers yields the desired result.

(e) Consider an arbitrary element $z\in m$.  For $q \gg 0$ and any
$q'$ we may replace $q$ by $qq'$, and so we know that
$c'z^{Nqq'} \in (y^{qq'}) + (I^{[qq']} \col R x^{qq'})$.

Multiply by $d^q$ to obtain
$$
d^qc'z^{Nqq'}\inc (d^qy^{qq'}) + d^q((I^{[qq']}
\Col x^{qq'})\inc (d^qy^{qq'}) +
(I^{[q']} \Col x^{q'})^{[q]}.
$$
where we are using (b) to show that, by (e),
$$d^q(I^{[qq']}\Col x^{qq'})  \inc
(I^{[q']}\Col x^{q'})\q.$$
Since $c'(dz^{Nq'})^q\inc (dy^{q'})^q + (I^{[q']} \Col x^{q'})^{[q]}$
for all $q'$ and all sufficiently large $q$, we have that
$dz^{Nq'}\in ((dy^{q'}) + (I^{[q']} \Col x^{q'}))^*$ for all $q'$.

Because $J$ is a strong test ideal, we have that
$$
Jdz^{Nq'}\in J(dy^{q'}) + J(I^{[q']} \Col x^{q'}).
$$
As in the final paragraph of 4.2, this yields that
$P(dz^{Nq'}, dy^{q'})\in (I^{[q']} \Col x^{q'})$
where $P(U,V) = U^k+r_1U^{k-1}V\,+\,\cdots\,+\,r_nV^k$ is a
homogeneous polynomial
of degree
$k$ in two variables over $R$ monic in $U$.
Factoring out $d^k = c$ we obtain that
$cP(z^{Nq'},y^{q'})\in (I^{[q']} \Col x^{q'})$,
and hence that $P(z^{Nq'},\,y^{q'})\in \Ap$ for all $q'$. It follows that
$z^{Nkq'}\in (y^{q'}) + \Ap$ for all $q'$.  Set $L$ equal to
$k$ times the minimal number of generators of $m^{Nk}$. Then
$m^{NLq'}\subseteq
(m^{Nk})^{[q']}$, and so $m^{NLq'}\subseteq (y^{q'}) + \Ap$ for all $q'$.\qed
\enddemo

\proclaim{Lemma 4.4 (push-up lemma)} Let $(R,\,m)$ be an excellent local
domain of positive prime characteristic $p$, let $I \inc R$ and let $x \in I$.
Let $y$ be an element of $R$ that is not in any associated prime, with
the possible exception of $m$, of
either  $(I+xR)\q$ or $I\q$ for any $q$. Let $P$ be a prime ideal in $T_I(x)$
of least height $h$, and let $J = I+yR$. If $h \leq \dim R - 2$,
then every minimal prime
$Q$ of $P+yR$ is in $T_J(x)$. \endproclaim

\demb{Proof} Let $c$ be a square locally stable test element.
It suffices to prove that $Q$ contains
$J\q \Col cx^q$ for all large $q$, since if so, then as
$Q$ is minimal over $P + Ry$, it will have to be minimal over
$J\q \Col cx^q$ as well.  To see this, note that if $Q'\subset Q$
strictly is a prime containing
$J\q \Col cx^q$, then $Q'$ will
contain some minimal prime over $\A$, and then $Q'$ will contain some prime
$Q''$ in $T_I(x)$. But the height of $Q''$ is at least $h$,
and the height of $Q'$ is at most $h$, so that $Q' = Q''$
is forced. But $y\in Q'$ then gives a contradiction as $y$ is
not contained in any associated prime of $I^{[q]}$ except possibly $m$.
Assume that $Q$ does not contain $((I,y)^{[q]} \Col cx^q)$. Then in $R_Q$,
$cx^q\in J^{[q]}_Q$, and we can write $cx^q - r_qy^q\in I^{[q]}_Q$ for some
$r_q\in R_Q$. Once we have localized at $Q$, which is properly
contained in $m$,  $y$ is not a zero divisor on $(I+xR)^{[q]}$,
since it is outside the  associated primes of
$(I+xR)\q$ other than $m$.
We therefore obtain that $r_q\in (I^{[q]} + x^q)_Q$. We then have that
for some $s_q \in R_Q$ there is an equation
$$
x^q(c - s_qy^q)\in I^{[q]}_Q.
$$
This implies that $x\in (I_Q)^*$.  (It suffices to work modulo a minimal
prime of the completion to check this, and we may apply Theorem
3.1 of [HH3].  The point is that the
$q\,$th root of the image of $c - s_qy^q$ will have
order approaching 0.) But then $x\in (I_{P})^*$.
Since $P\in T_I(x)$ this is a
contradiction, which proves that  $Q\in T_J(x)$. \qed\enddemo

We are now ready for the main result of this section.

\proclaim{Theorem 4.5} Let $(R,m)$ be a local domain with countable
prime avoidance, module-finite, and generically \'etale
over an excellent regular local ring
$A$. Suppose that $R$ has a strong test ideal $J \not= 0$ with $k$
generators.
Suppose also that for all $P\ne m$, $R_P$ has the property that the tight
closure of ideals commutes with localization. Furthermore assume that $R$
satisfies (LC). Fix a square locally stable test element $c$ that is
the $k\,$th power of an element of $A\0$ multiplying $R$ into
$A[R^p]$, as in Lemma 4.3(e).
Then the following are equivalent:
{\roster
\item Tight closure of ideals in
$R$ commutes with localization.

\item For all ideals $I$ and all $x\in R$
there exists an $\epsilon > 0$ such that for all $Q\in T_I(x)$,
$$ \ell(q,I,x,Q) = \lambda((R_Q/(\A)_Q)) > \e q^{\text{dim}(R_Q)}$$
for all $q \gg 0$.
\endroster}
\endproclaim

\demb{Proof} We first prove $(1) \imp (2)$. Assuming
(1), Theorem 3.5 gives us that $T_I(x)$ is a finite set.
Hence it suffices to prove that for a fixed
$Q\in T_I(x)$, $\text{lim inf}_q\lambda_q(Q) > 0$.

Let $\mu$ denote the number of generators of $R$ as an $A$-module.
Set $P = Q\cap A$.  Let $W = A-P$.  Let  $s$ be the degree of
the residue field extension $[R_Q/QR_Q:A_P/PA_P]$:  note that
for $R_Q$-modules length over $A_P$ is $s$ times the length
over $R_Q$.

Since $x \notin (I_Q)^*$, we can
choose $q'$ such that $cx^{q'} \notin I^{[q']}_Q$ and then we have:\smallskip

\noindent $s\mu \cdot \ell(q,I,x,Q) \geq
   \mu \cdot \lambda_{A_P}(A_P/((\A)_Q\cap A_P) \geq
     \lambda_{A_P}(R_W/((\A)_Q\cap A_P)R_W) \geq$

\noindent $\lambda_{R_Q}(R_Q/((\A)_Q\cap A_P)R_Q)\geq
    \lambda_{R_Q}(R_Q/((\Ap)_Q)^{[q/q']})$

\medskip
\noindent by Lemma 4.3(c). But then
$\ell(q,I,x,Q) \geq Cq^{\text{dim}(R_Q)}$ for some $C >0$
for large $q$ since $\lambda_{R_Q}(R_Q/((\Ap)_Q)^{[q/q']})$ is
asymptotic to a positive constant
times $(q/q')^{\text{dim}(R_Q)}$ by the standard theory of
Hilbert-Kunz functions
[Mon]. This completes the direction `$\imp$'.

We now consider the much more difficult converse direction.
Assume the conditions
in (2).  Notice that because we have assumed that localization
holds for any proper localization of $R$, we know that
condition C2 holds for all $I$, $x$, $c$ (if $Q = m$ we know the
condition (LC), which is stronger).  Suppose that the result is false, and
choose $I$ maximal such that the
tight closure of $I$ does not
commute with localization. By Theorem 3.5 and Remark 3.6,
there exists an element $x$
(and we may assume that $x \in (I_Q)^* - (I^*)_Q$ for some prime $Q$)
such that either C1 or C2 fails. Since we know the latter,
there must
be such an $x$ with $T_I(x)$ infinite.  Notice that if we replace
$x$ by $zx$ for any choice of $z \notin Q$,  the element $zx$ is
still in $(I_Q)^* - (I^*)_Q$, and so $T_I(zx)$ is still infinite.
We shall make several such replacements that will force increasingly
controlled behavior on $T_I(zx)$.  Each time, we change notation
and write $x$ for what is really $zx$.  We shall eventually obtain
a contradiction.

First apply condition (LC) to choose $N$ such that $m^{Nq}$
kills $H^0_m(R/I^{[q]})$ for all $q$, and choose $z\in m^N$ not
in $Q$,  and not in any associated prime, except $m$, of any
of the countably many ideals $I\q$.
Then the remark above shows that
$T(zx)$ is still infinite.  Moreover,  $m$ is not associated to
$I^{[q]} \Col uz^q$ for
all $q$ and any choice of $u$ in $R$, by 3.4(g).

We replace $x$ by $zx$ and can assume from now on  that $m$ is not
associated to
$\A$ for all $q$. Furthermore if we again replace $x$ by $z'x$ for some
$z'$, this remains true.

We next replace $x$ by $zx$ for $z \notin Q$ so as to maximize the
least height $h$ of a prime in $T_I(xz)$.  Thus, without loss of
generality we may assume that for all $z \in R-Q$,  replacing
$x$ by $zx$ does not clear all the primes of height $h$. Notice
that by 3.4(f) as we make such replacements the least height
occurring can only increase.

Since C2 holds it is immediate from the clearing lemma and the
remark above that there must be infinitely many primes of
height $h$.

We shall next show that $h = d-1$.  Assume $h \leq d-2$.
We can choose an element $y$ of $m$ not in any associated prime, except
possibly $m$, of any of the ideals $I\q$, $(I+Rx)\q$, and then
the push-up lemma shows that for every prime $P$ of height $h$
in $T_I(x)$, every minimal prime $P'$ of $P+yR$ is in $T_{(I+Ry)}(x)$.
But localization holds in $R_{P'}$, so that only finitely many primes
in $T_I(x)$
can lie inside $P'$ (note that
$T_I(x)\cap \text{Spec}(R_{P'}) = T_{I_{P'}}(x)$:  see 3.3(c)).
On the other hand, the maximality of $I$ forces $T_{(I,y)}(x)$ to be finite
also. As we have shown that every one of the infinitely many primes
in $T_I(x)$ of height $h$ lies inside some $P\in T_{(I,y)}(x)$, this
contradiction proves that $h = d-1$.

Henceforth, we may assume that $h = d-1$. We choose $y$ sufficiently
general as above.  Recall that $c$ satisfies the hypothesis
of Lemma 4.3(e):
we shall need this below.  We write $J\sat$ for the inverse image
of $H^0_m(R/J)$ in $R$, i.e., for the union of all the ideals
$J\Col m^N$ as $N$ varies.

By Noetherian induction, the ideal $(I,y)$ satisfies conditions
C2$^*$ and C1 and we
know that
the ideals $(I,y)^{[q]} \Col cx^q$ are either the whole
ring or are $m$-primary
(since $y$ is not in any of the primes of $T_I(x)$),
and there exists an integer $N$ such that for all $q$
$m^{Nq}\subseteq (I,y)^{[q]} \Col cx^q$. Choose a general element
$z\in m^{Nq}$
and write $czx^q - r_qy^q\in \I$ for some $r_q\in R$. Then $r_q\in
(\I + (x^q)) \Col y^q\inc (\I + (x^q))\sat$ by the choice of $y$. Hence the
quotient $((\I + (x^q)) \Col y^q)/(\I + (x^q))$ is contained in
$H^0_m(R/(\I + (x^q)))$
and there exists an integer $B$ such that $m^{Bq}$ annihilates this quotient.
In particular, $m^{Bq}r_q\in (\I + (x^q))$, and for
general $u\in m^{Bq}$ we can write
$$czux^q - s_qx^qy^q\in \I$$ for some $s_q\in R$. Then
$x^q(czu - s_qy^q)\in \I$ and it follows that
$$
cm^{(N+B)q}\inc (\I \Col x^q) + (y^q)
$$
for all $q$.

We now apply Lemma 4.3(e) to conclude that there is a constant
$D$ such that for all $q$,
$$
m^{Dq}\inc (y^q) + \A.
$$
In particular, $\lambda(R/((y^q) + (\A)))\leq \lambda(R/m^{Dq})$ and for
large $q$ it follows that there is a constant $C > 0$ such that
$$
\lambda(R/((y^q) + (\A))\leq Cq^d.
$$

Fix $\e > 0$ as in the statement of the theorem.
Let $N_q$ be the number of
minimal primes above $\A$.
We claim that
$$N_q\leq C/\e$$
where $C$ is as in the above paragraph.  Since every minimal prime
above $\A$ is also minimal above $\Ap$ for
$q'\geq q$ in our case, this means that the minimal primes above
$\A$ stabilize for large $q$, and hence implies that $T_I(x)$ is a
finite set.

Recall that $m$ is not associated to $\A$. The associativity
formula for multiplicities gives that
$$ \lambda(R/((y^q) + \A) =  e(y^q\sc R/(\A)) = qe(y\sc R/(\A))$$
and
$$
qe(y\sc R/(\A)) =
q(\sum_{P\in \text{Min}(\A)} e(y \sc R/P)\lambda(R_P/(\A)_P)
$$
and hence
$$
Cq^{d-1} \geq \sum_{P\in \text{Min}(\A)}
\lambda(R_P/(\A)_P) = \sum_{P\in \text{Min}(\A)} l(q,I,x,P)).
$$
Thus $Cq^{d-1}\geq N_q\cdot \e q^{d-1}$
and hence $N_q \leq C/\e$ as claimed.\qed
\enddemo
\vfill\eject

\Bigskip
{\chbf
\centerline{\hdbf 5. FURTHER REMARKS AND QUESTIONS}
}
\bigskip

\demb{Discussion 5.1: finding specific test exponents}
Tight closure is known to commute with
localization in many specific cases:  under mild conditions on
the ring this is true for ideals generated by monomials in
parameters and ideals $I$ such that $R/I$ has finite phantom
projective dimension. We refer the reader to [AHH, \S8] for
a detailed discussion of various results.  The known results on
when tight closure
commutes with localization therefore imply the existence of
test exponents for many ideals.  However, little is known about
how to determine a specific test exponent for a given ideal $I$.
We want to raise this as a problem.  If one has a specific test
exponent for $c, I$ then to test whether $u \in I^*$ one need
only test whether $cu^q \in I\q$ for that one value of $q$. We
believe that the best hope for giving a useful algorithm
for testing when an element is in the tight closure of an ideal
lies in this direction.  It would be of considerable interest
to solve the problem of determining test exponents effectively
even for parameter ideals.
\enddemo

\demb{Discussion 5.2: algorithmic testing for tight closure} We
want to point out that in certain instances there is an
algorithm, in a technical sense, for testing whether specific
elements are in a tight closure.  We do not believe that this
particular method will ever be implemented.  In any given instance
where it may be applied, it does eventually terminate, showing
that the specific element $x$ either is or is not in the
tight closure of $I$.  However, we do not have a way of estimating
{\it a priori} how long testing may need to go on before the
algorithm terminates.

The method may be applied to ideals $I$ such that the tight closure
of $I$ is the same as the plus closure of $I$.  We review the latter
notion.  Suppose that $R$ is a domain.   If $I \inc J \inc R$ and there
is an integral extension (equivalently,
a module-finite extension) $S$
of $R$ with  $J \inc IS$, then
$J \inc I^*$. If $R\pl$ denotes
the integral closure of $R$ in
an algebraic closure of
its fraction field (the {\it
absolute integral closure}:  see [HH4] and [Smith] for
further discussion), we can
let  $I^+ = IR\pl \cap R$ and
then  $I \inc I^+ \inc I^*$. Is $I^* = I^+\,$?
Except in trivial cases where $I = I^*$
for all $I$, we do not know whether this is true in
any normal domain.
By a hard theorem (cf.\ [Smith]), it is true for
parameter ideals, and a result of
Aberbach [Ab]  permits one to extend this
to ideals $I$ such that $R/I$ has
{\it finite phantom projective dimension}.
The point we want to make is that for ideals such that
$I^* = I^+$, which may well be all ideals,  there is an
algorithm of sorts.

Let $R$ be a countable Noetherian
domain of prime characteristic $p > 0$ in which basic
operations can be performed algorithmically,
with a known test element $c$, and such that one
can test algorithmically for membership
in an ideal in polynomial rings over $R$,
e.g., a finitely generated domain over a finitely
generated field. \medskip

\noindent {\bf Fact.} {\sl If $R$ is as above and
$I \inc R$ satisfies $I^* = I^+$, then one
can test algorithmically whether $y \in R$
is in $I^*$. In particular, one has such a test
if $R$ is an affine domain over a finitely  generated field
and $I$ is generated by monomials in elements
$\vect z d$ generating an ideal of height $d$.}
\medskip

Here is the idea of the algorithm:  one can effectively enumerate all
the algebras
$\vect S n,\,\ldots$  that are
module-finite over $R$.  Alternately test
whether $y \in IS_n$ and whether $cy^{p^n} \in I^{[p^n]}$.  If
$y \in I^* = I^+$ the former test eventually succeeds, while
if $y \notin I^*$, the latter test eventually
fails. \qed

Of course, this method is awful:  this algorithm only gives
emphasis to the problem of effective determination of test
exponents.
\enddemo

\demb{Discussion 5.3: uniform test exponents} Let $R$ be
reduced and finitely generated over an excellent local ring.  So far as
we know it is possible that for a given locally stable test element
$c$ there exists a test exponent valid for all ideals $I$ simultaneously.
It would suffice to give such an exponent for all $m$-primary
tightly closed ideals as $m$ varies, and even for those which
are maximal with respect to the property of being tightly closed
and not containing a given element of the ring, since every
tightly closed ideal is an intersection of such ideals.

A more modest question that seems more approachable is
whether, given $c$,  there exists a single test exponent for
all ideals containing a given $m$-primary ideal $J$, because
then one can construct a moduli space for the set of ideals.
\enddemo

\demb{Discussion 5.4} It is reasonable to ask whether
localization can be proved for suitable local domains, e.g., those
that are excellent, have countable prime
avoidance, and are generically \'etale over a regular local ring,
if for all ideals $I \inc R$ and all $x \in R$,  there exist constants
$\epsilon$, $\delta$
greater than 0 such that for all $q$ and for all $Q\in T_I(x)$,
$\delta q^{\dim(R_Q)} \geq \ell(q,I,x,Q) \geq \epsilon q^{\dim(R_Q)}$.
If localization holds the first inequality can be deduced from C2, while
the second was shown in the proof of Theorem 4.5 to follow from C1.
Thus, these conditions are necessary.  It also  seems reasonable to
ask  whether $\ell(q,I,x,Q)/q^{\dim(R_Q)}$ approaches a
(necessarily positive) limit as $q \to \infty$.
\enddemo
} 

\Bigskip
\Bigskip

{
\smc
\baselineskip = 10 pt
\settabs 7 \columns
\quad\bigskip
\+Department of Mathematics           &&&&Department of Mathematics\cr
\+University of Michigan              &&&&University of Kansas\cr
\+Ann Arbor, MI 48109--1109           &&&&Lawrence, KS 66045\cr
\+USA                                 &&&&USA\cr
\smallskip
\+E-mail:                             &&&&E-mail:\cr
\vskip 1.5 pt plus .5 pt minus .5 pt
{\rm
\+hochster\@math.lsa.umich.edu        &&&&huneke\@math.ukans.edu\cr
}
}

\enddocument